\documentclass[12pt,a4paper,twoside,reqno]{amsart}

\usepackage{amssymb,amsmath,amstext,amsthm,amsfonts}
\usepackage[dvips]{graphicx}
\usepackage[ansinew]{inputenc} 
\usepackage{fullpage}

\newtheorem{maintheorem}{Theorem}

\newcommand{\cmt}{\begin{maintheorem}}
\newcommand{\fmt}{\end{maintheorem}}

\newtheorem{T}{Theorem}[section]
\newcommand{\cte}{\begin{T}}
\newcommand{\fte}{\end{T}}

\newtheorem{Corollary}[T]{Corollary}
\newcommand{\cco}{\begin{Corollary}}
\newcommand{\fco}{\end{Corollary}}

\newtheorem{Proposition}[T]{Proposition}
\newcommand{\cpr}{\begin{Proposition}}
\newcommand{\fpr}{\end{Proposition}}

\newtheorem{Lemma}[T]{Lemma}
\newcommand{\cle}{\begin{Lemma}}
\newcommand{\fle}{\end{Lemma}}

\theoremstyle{remark}
\newtheorem{Example}[T]{Example}
\newcommand{\cex}{\begin{Example}}
\newcommand{\fex}{\end{Example}}

\newtheorem{Remark}[T]{Remark}
\newcommand{\cre}{\begin{Remark}}
\newcommand{\fre}{\end{Remark}}

\newtheorem{Definition}[T]{Definition}
\newcommand{\cde}{\begin{Definition}}
\newcommand{\fde}{\end{Definition}}

\newcommand{\Chi}{{\bf 1}}

\newcommand{\mcup}{\mbox{$\bigcup$}}

\newcommand{\be} {\beta}        
       
\newcommand{\de} {\delta}       
\newcommand{\ep} {\epsilon}
\newcommand{\e} {\epsilon}
\newcommand{\vare}{\varepsilon}

\renewcommand{\th} {\theta}

\newcommand{\si} {\sigma}

\def \RR {{\mathbb R}}
\def \ZZ {{\mathbb Z}}
\def \NN {{\mathbb N}}

\newcommand{\dem}{\begin{proof}}
\newcommand{\cqd}{\end{proof}}

\newcommand{\Leb}{m}

\newcommand{\dist}{\operatorname{dist}}

\newcommand{\leb}{m}

\newcommand{\qand}{\quad\text{and}\quad}
\newcommand{\underi}{{\underline{i}}}

\newcommand{\cc}{{\mathcal C}}

\newcommand{\cf}{{\mathcal F}}

\newcommand{\cn}{{\mathcal N}}

\newcommand{\cp}{{\mathcal P}}

\title[Strong statistical stability]{Strong statistical stability\\
of non-uniformly expanding maps}

\author{Jos\'e F. Alves}

\address{CMUP\\
Rua do Campo Alegre 687, 4169-007 Porto, Portugal}
\email{jfalves@fc.up.pt}
\urladdr{http://www.fc.up.pt/cmup/jfalves}
\date{\today}

\thanks{Work partially supported by FCT through CMUP}



\begin{document}

\begin{abstract}
We consider families of transformations in multidimensional Riemannian manifolds with non-uniformly expanding
behavior. We give sufficient conditions for the continuous variation (in the $L^1$-norm) of the densities of
absolutely continuous (with respect to the Lebesgue measure) invariant probability
 measures for those transformations.
\end{abstract}

\maketitle


\section{Introduction} \label{se.introduction}

In this work we address ourselves to the study of the statistical stability of certain classes of chaotic
dynamical systems. We are particularly interested in the statistical stability of systems  displaying
non-uniformly expanding behavior on the growth of the derivative for most of its orbits.

To be more specific, let  $f : M \rightarrow M$ be some discrete-time dynamical system of a compact Riemannian
manifold $M$, and let $m$ be a volume form that we call \emph{Lebesgue measure}.  {\em Sinai-Ruelle-Bowen (SRB)
measures} or {\em physical measures} are probability measures that characterize asymptotically, in time average,
a large set of orbits of the phase space; these are defined precisely in \eqref{eq.SRB} below. It is a difficult
problem to verify the existence of these measures for general dynamical systems.

 By the {\em statistical
stability} of a system, we mean continuous variation of the SRB measures under small modifications of the law
that governs the system. Using Birkhoff's Ergodic Theorem, one possible way for finding SRB measures for a map
$f$ is by proving the
existence of ergodic absolutely continuous 
$f$-invariant probability measures.

Systems displaying uniformly expanding behavior  have been exhaustively studied in the last decades, and several
results on the existence of SRB measures and their statistical stability have been obtained, starting with
Sinai, Ruelle and Bowen; see \cite{Si,R, BR,Bow75} and also \cite{R2,Ki1, Ki2, Yo}.

 The existence of SRB measures for many
one-dimensional maps with non-uniformly expanding behavior has been established in the pioneer work of Jakobson
\cite{J}; see also \cite{BC1,BY,BaV}. Viana introduced in \cite{V}  an open class of transformations in higher
dimensions with non-uniformly expanding behavior for most of its orbits. The existence of SRB measures for Viana
maps has been proved in \cite{Al}. Motivated by the results in \cite{V} and~\cite{Al}, general conclusions on
the existence of SRB measures for non-uniformly expanding dynamical systems   are drawn in \cite{ABV}.

The statistical stability of the systems introduced in \cite{V} has been proved in \cite{AV}, in a strong sense:
convergence of the densities of the SRB measures in the $L^1$ norm. The proof uses in an important way
geometrical features of the system, and could not be immediately extended to more general classes of
non-uniformly expanding maps. Some results in this direction were obtained in \cite{Al2}, but in a weak sense:
convergence of the measures in the weak* topology.

In this work we give sufficient conditions for the strong  statistical stability of certain classes of
non-uniformly expanding maps. These conditions are naturally verified by the maps introduced in \cite{V}, as
shown in \cite{AV}, and by a class of non-uniformly expanding local diffeomorphisms introduced in \cite{ABV}
that we include at the end of this work.


\subsection{Non-uniformly expanding maps}\label{se.NUE}

Let $f\colon M\to M$   be a continuous map
which is local diffeomorphism in the whole manifold  except  in  a
set of critical  points $ \cc \subset M$.

 \cde\label{d.nd}  We say that $ \cc $ is
{\em non-degenerate}\index{non-degenerate!critical set} if the
following conditions hold. The first one
 says that
 $f$  {\em behaves like a power of the distance}
 to $ \cc $:
 there are $B>1$ and $\be>0$ such that for every $x\in
 M\setminus\cc$
\begin{enumerate}
 \item[(s$_1$)]
\quad $\displaystyle{{B}^{-1}\dist(x,\cc)^{\be}\leq 
{\|Df(x)v\|}
\leq B\dist(x,\cc)^{-\be}}$, for all $v\in T_x M$ with $\|v\|=1$.
\end{enumerate}
Moreover, we assume that  $  \log|\det Df| $ and $ \log
\|Df^{-1}\| $ are \emph{locally Lipschitz} in $ M \setminus \cc$,
with Lipschitz constant depending on the distance to $\cc$: for
every $x,y\in M\setminus \cc$ with $\dist(x,y)<\dist(x,\cc)/2$ we
have
\begin{enumerate}
\item[(s$_2$)] \quad $\displaystyle{\left|\log\|Df(x)^{-1}\|-
\log\|Df(y)^{-1}\|\:\right|\leq
\frac{B}{\dist(x,\cc)^{\be}}\dist(x,y)}$;
 \item[(s$_3$)]
\quad $\displaystyle{\left|\log|\det Df(x)|- \log|\det
Df(y)|\:\right|\leq \frac{B}{\dist(x,\cc)^{\be}}\dist(x,y)}$.
 \end{enumerate}
 \fde

Given $\delta>0$ and $x\in M\setminus\cc$ we define the {\em
$\delta$-truncated distance\/} $\dist_\delta(x,\cc)=\dist(x,\cc)$,
if  $\dist(x,\cc)<\delta$, and $\dist_\delta(x,\cc)=1$, otherwise.
%

\cde\label{def.NUE} Let $f\colon M\to M$ be a local diffeomorphism
outside a non-degenerate critical set $\cc$. We say that $f$ is
{\em non-uniformly expanding}  if:
\begin{itemize}
\item there is $\lambda>0$ such that for every $x\in M$
 \begin{equation}\label{NUE}
    \limsup_{n\to\infty}\frac{1}{n}\sum_{i=0}^{n-1}
    \log \|Df({f^{i}(x))}^{-1}\|<-\lambda;
\end{equation}
\item  for every $\epsilon>0$ there exists $\delta>0$ such that
for every $x\in M$
\begin{equation} \label{e.faraway1}
    \limsup_{n\to+\infty}
\frac{1}{n} \sum_{j=0}^{n-1}-\log \dist_\delta(f^j(x),\cc)
\le\epsilon.
\end{equation}
\end{itemize}
We will often refer to \eqref{e.faraway1} by saying that orbits
have {\em slow recurrence} to the critical set $\cc$. When
$\cc=\emptyset$ we simply ignore the slow recurrence condition.
 \fde

\cre\label{remark} Slow recurrence condition  is not needed in all its strength. In fact, the only place where
we will be using \eqref{e.faraway1} is in the proof of Proposition~\ref{pr.hyperbolic1}. As we shall see, it is
enough that~\eqref{e.faraway1} holds for some sufficiently small $\ep>0$ and conveniently chosen $\delta>0$; see
Remark~~\ref{r.strong}. \fre

A Borel probability measure $\mu$ on the Borel sets of $M$ is said to be an {\em SRB measure\/} if there exists
a positive Lebesgue measure set of points $z\in M$ for which
\begin{equation}
\label{eq.SRB} \lim_{n\to+\infty}\frac{1}{n}\sum_{j=0}^{n-1}
\varphi(f^j(z)) =\int \varphi\,d\mu
\end{equation}
for any continuous function $\varphi:M\to\RR$. The set of points $z\in M$ for which this holds is called the
\emph{basin} of $\mu$. It was proved in \cite{ABV} that non-uniformly expanding maps possess SRB measures.

If $f\colon M\to M$ is non-uniformly expanding, then by
\eqref{NUE} the \emph{expansion time}\index{time!expansion}
function
\begin{equation}\label{de.exptime}
\mathcal E(x) = \min\left\{N\ge1: \frac1n\sum_{i=0}^{n-1} \log
\|Df({f^{i}(x))^{-1}}\| \leq -\lambda, \quad\text{for all $n\geq
N$}\right\}
\end{equation}
is defined and finite almost everywhere in $ M $.  Then, according
to Remark~\ref{remark}, we fix $\vare>0$ and $\delta>0$ as in
\eqref{e.faraway1}. The \emph{recurrence
time}\index{time!recurrence} function
\begin{equation}\label{de.slowrec}
\mathcal R(x) = \min\left\{N\ge 1: \frac1n\sum_{i=0}^{n-1} -\log
\dist_\delta(f^j(x),\cc) \leq \varepsilon , \quad \text{for all
}n\geq N\right\}
\end{equation}
is also defined and finite almost everywhere in $ M $.  We define
the {\em tail set}
\begin{equation}\label{tailset}
    \Gamma_n=\big\{x: \mathcal E(x) > n \ \text{ or } \ \mathcal R(x) > n \big\}.
\end{equation}
This is the set of points which at time $ n $ have not yet
achieved either the uniform exponential growth of derivative or
the uniform slow recurrence. If $\cc=\emptyset$, we ignore the
recurrence time function
 in the definition of $\Gamma_n$.

\subsection{Statistical stability}
Let $\cf$ be a family of $C^k$ maps ($k\ge 2$) from a $d$-dimensional manifold $M$ into itself,
and endow $\cf$ with the $C^k$ topology. We assume that each $f\in\cf $  admits a unique absolutely
continuous $f$-invariant probability measure $\mu_f$ in $M$. 

\cde\label{de.sstable} We say that $f_0\in\cf$ is {\em (strongly) statistically
stable}\index{stable!statistically}, if
 $
 \cf\ni f \mapsto {d\mu_f}/{dm}
 $
 is continuous at $f_0$, with respect to the $L^1$-norm on the space of densities. \fde

We  assume that the maps in a neighborhood of $f_0$ satisfy the
following non-degeneracy condition: given any $\e>0$ there exists
$\de>0$ such that
\begin{equation}
\label{eq.oito} m(E) \le \de \quad\Rightarrow \quad m(f^{-1}(E))
\le \e
\end{equation}
for any measurable subset $E\subset M$ and any $f\in\cf$. This can often be enforced by requiring some jet of
order $l\le k$ of $f_0$ to be everywhere non-degenerate. This is obviously satisfied whenever we consider local
diffeomorphisms.

\cde\label{d.uniform} We say that $\cf$ as above is a {\em uniform family}\index{uniform family} if the
$B,\beta$ as in Definition~\ref{d.nd}, and $\vare, \delta, \lambda$ as in Definition~\ref{def.NUE} (cf.
Remark~\ref{r.strong}) can be chosen uniformly in  $\cf$. \fde

\cmt\label{te.uniform} Let $\cf$ be a uniform family of $C^k$
($k\ge 2$) non-uniformly maps for which non-degeneracy condition
\eqref{eq.oito} holds.  Assume that there are $C>0$ and $\gamma>1$
 such that
 $m(\Gamma^f_n)\le Cn^{-\gamma},$ for all $n\ge 1$ and $f\in\cf$.
 Then every $f\in \cf$ is statistically stable.
 \fmt
Condition \eqref{eq.oito} is needed just because we are going to use \cite[Theorem A]{AV}; see
Theorem~\ref{t.A2} below.


\section{Piecewise expanding induced maps}\label{se.piecewise}

One possible way for proving the existence of invariant measures for certain dynamical systems may be  by
choosing conveniently some region in the phase space and studying an induced return map to that region. This
method can also be efficient in proving the absolute continuity of those measures. In this section we are
particulary interested in the study of those return maps.

\subsection{Markovian return maps}\label{se.return}

Let $f$ be a map from a  Riemannian manifold $M$ into itself, and let $ F:\Delta\rightarrow \Delta$ be a {\em
return map}\index{return map} for $f$ in some topological disk in $\Delta\subset M$. This means that there is a
countable partition $\cp$ of a full Lebesgue measure subset of $\Delta$, and there exists a {\em return
time}\index{time!return} function $R\colon\cp\rightarrow\ZZ^+$ such that $  F\vert_U=f^{R(U)}\vert_U $ for each
$ U\in \cp. $

\cde\label{de.expanding} We say that $F$ is a  {\em piecewise
expanding Markovian map\/} if there is a countable partition $\cp$
into open sets of a full Lebesgue measure subset of $\Delta$ such
that:
 \begin{enumerate}
 \item {\em Expansion:}\index{expansion!uniform}
 there is $0<\kappa<1$ such that for each  $U\in\cp$ and $x\in U$
 $$\|D F(x)^{-1}\| <\kappa.$$
 \item {\em Bounded distortion:\/}
there is  $K>0$ such that for each $U\in\cp$ and $x,y\in U$
 \[
\log\left|\frac{\det DF(x)}{\det DF(y)}\right| \leq K \dist(F(x),
F(y)).
\]
  \item {\em {Markov}:}
$F\vert_U$ is a $ C^{2} $ diffeomorphism  onto $\Delta$, for each
$U\in\cp$.
\end{enumerate}
\fde

If $F\colon\Delta\to\Delta$ is a $C^2$ piecewise expanding Markovian map, then it has some absolutely continuous
invariant measure $\mu_F$. Moreover, the density of $\mu_F$  is uniformly bounded by some constant; see e.g.
\cite[Theorem 1]{Y2}. Defining
 \begin{equation}\label{cau}
 \mu^{\ast}_f
=\sum_{j=0}^{\infty}f_{\ast}^j\left(\mu_F\mid \{R>j\}\right),
 \end{equation}
it is straightforward to check that $\mu^{\ast}_f$ is an  absolutely continuous $f$-invariant
measure, which is finite  whenever $R\in L^1(\Delta)$. 

\subsection{Statistical stability}

Let $\cf$ be a family of $C^k$ maps ($k\ge 2$) from the manifold $M$ into itself, and assume that we may
associate to each $f\in\cf$ a piecewise expanding return  map $ F_f\colon \Delta\rightarrow \Delta$  as in
Definition~\ref{de.expanding}. For each $f\in\cf$, let $\cp_f$ denote the partition into domains of smoothness
of $F_f$ and $R_f\colon\cp_f\rightarrow \ZZ^+$ be the corresponding return time. We assume that $R_f\in
L^1(\Delta)$ for each $f\in\cf$, which then implies that if $\mu_F$ is the absolutely continuous $F_f$-invariant
probability measure, then
 $
 \mu_f^{\ast}
=\sum_{j=0}^{\infty}f_{\ast}^j\left(\mu_F\mid \{R_f>j\}\right)
 $ 
is an absolutely continuous $f$-invariant finite  measure. 
In our setting of Markovian maps, the statement of \cite[Theorem A]{AV} can be simplified.

\cte \label{t.A2} Let $\cf$ be as above, and suppose that every $f\in\cf$ admits a unique absolutely continuous
invariant probability measure $\mu_f$. Suppose that each $f_0\in\cf$ satisfies:
 \begin{enumerate}
 \item[(u$_1$)]
 Given $\epsilon>0$ there is  $\de>0$ such that for any $f\in\cf$
 $$
\|f-f_0\|_{C^k}<\de\quad\Rightarrow\quad\|R_f-R_{f_0}\|_1<\epsilon.
 $$
\item[(u$_2$)] $\kappa$, $K$ as in Definition \ref{de.expanding} may be taken uniformly in a neighborhood of
$f_0$ in $\cf$.
 \end{enumerate}
 Then $f_0$ is statistically stable. \fte

\cre The bounded distortion condition used in \cite[Theorem A]{AV} is satisfied in our context, as we shall see
in Lemma~\ref{le.bd}. Moreover, the assumption on the constants $\beta$ and $\rho$ as in condition (U3) of
\cite{AV} is trivially satisfied. In the non-Markovian case treated in \cite[Theorem A]{AV}, one can only assure
that the density of $\mu_F$ belongs to $L^p(\Delta)$ for some $p>1$. This implies that convergence of $R_f$ to
$R_{f_0}$ has to be taken in the norm of $L^q(\Delta)$ with $p^{-1}+q^{-1}=1$. Since in our case the density
belongs to $L^\infty(\Delta)$ we may take the convergence of $R_f$ to $R_{f_0}$  in the $L^1$-norm, by a usual
H\"older inequality argument. \fre

Under the assumptions of the Theorem~\ref{t.A2},  the unique absolutely continuous invariant probability measure
is necessarily equal to the normalization of $\mu_f^*$, i.e.
 $\mu_f=\mu_f^*/{\mu_f^*(M)}.$
Thus for proving Theorem~\ref{te.uniform} we just have to show
that  conditions (u$_1$) and (u$_2$) hold for families $\cf$ as in
 Theorem~\ref{te.uniform}.


\section{Hyperbolic times and bounded distortion}

In this section we present some results on the existence of hyperbolic times for non-uniformly expanding maps
and  distortion properties at hyperbolic times. Although these results have essentially been all proved in
\cite{ABV}, we include some proofs here in order to see how the   constants depend on one another.

 \cde   Fix $B>1$ and $\beta>0$ as in Definition~\ref{d.nd},
and take $b>0$ such that $2b < \min\{1,\beta^{-1}\}$. Given $\sigma<1$ and $\delta>0$, we say that $n$ is a {\em
$(\sigma,\delta)$-hyperbolic time}\index{time!hyperbolic} for a point $x\in M$ if for all $1\le k \le n$,
 \begin{equation}\label{d.ht}
\prod_{j=n-k}^{n-1}\|Df(f^j(x))^{-1}\| \le \sigma^k \qand \dist_\delta(f^{n-k}(x), \cc)\ge \sigma^{b k}.
 \end{equation}
In the case $\cc=\emptyset$ the definition of $(\sigma,\delta)$-hyperbolic time reduces to the first condition
in \eqref{d.ht} and we simply call it a {\em $\sigma$-hyperbolic time}.
  \fde

\cle\label{leaux} Given $\delta>0$ fix $\delta_1=\delta_1(B, \beta, \sigma, \delta) > 0$ so that $4\delta_1 <
\delta$ and $4 B \delta_1 < \delta^\beta|\log\sigma|$. If $n$ is a $(\sigma,\delta)$-hyperbolic time for $x$,
then
$ 
\|Df(y)^{-1}\| \le \sigma^{-1/2} \|Df(f^{n-j}(x))^{-1}\|
$ 
for any $1 \le j < n$ and any point
$y$ in the ball of radius $2\delta_1 \sigma^{j/2}$ around $f^{n-j}(x)$. \fle

\dem Since $n$ is a $(\sigma,\delta)$-hyperbolic time for $x$ we have  $\dist_\delta(f^{n-j}(x),\cc) \ge
\sigma^{j}$ for any $1 \le j < n$.
 According to the definition of the truncated distance, this means that
$$
\dist(f^{n-j}(x),\cc)= \dist_\delta(f^{n-j}(x),\cc) \ge
\sigma^{bj} \quad\text{or else}\quad \dist(f^{n-j}(x),\cc) \ge
\delta.
$$
In either case, we have $\dist(y,f^{n-j}(x)) <
\dist(f^{n-j}(x),\cc)/2$ for any $1 \le j < n$, because we chose
$b<1/2$ and $\delta_1 < \delta/4 < 1/4$. Therefore, we may use
(s$_2$) to conclude that
$$
\log\frac{\|Df(y)^{-1}\|}{\|Df(f^{n-j}(x))^{-1}\|} \le B
\frac{\dist(y,f^{n-j}(x))}
           {\dist(f^{n-j}(x)),\cc)^\beta}
\le B \frac{2\delta_1 \sigma^{j/2}}
           {\min\{\sigma^{b\beta j},\delta^\beta\}}.
$$
Since $\delta$ and $\sigma$ are smaller than $1$, and we took
$b\beta<1/2$, the term on the right hand side is bounded by
$2B\delta_1 \delta^{-\beta}$. Moreover, our second condition on
$\delta_1$ means that this last expression is smaller than $\log
\sigma^{-1/2}$.
\cqd

%
%

  \cpr \label{pr.distorcao} Let $0<\sigma<1$ and $\delta>0$. 
If $n$ is a $(\sigma,\delta)$-hyperbolic time for $x$, then there exists a neighborhood $V_n$ of $x$ such that:
\begin{enumerate}
\item $f^n$ maps $V_n$ diffeomorphically onto the ball of radius $\delta_1$ around $f^n(x)$;
 \item for each $x\in V_n$  we have
 $\|Df^n(x)^{-1}\|\le \sigma^{n/2}$;
 \item for all $1\le k <n$ and $y, z\in
V_n$, $$ \dist(f^{n-k}(y),f^{n-k}(z)) \le \sigma^{k/2}\dist(f^{n}(y),f^{n}(z)). $$
\end{enumerate}
  \fpr

\dem See \cite[Lemma 5.2]{ABV}. \cqd

We shall  refer to the sets $ V_n$ as \emph{hyperbolic pre-balls} and to their images $ f^{n}(V_n) $ as
\emph{hyperbolic balls}. Notice that the latter are indeed balls of radius $ \delta_1>0 $.

\cle \label{l.pliss} Given $0<c_1<c_2<A$ let $\theta=(c_2-c_1)/(A-c_1)$. Take $a_1\le A, \ldots, a_N\le A$ such
that $ \sum_{j=1}^{N} a_j \ge c_2 N .$ Then there are $l > \theta N$ and $1 < n_1 < \cdots < n_l \le N$ so that
$
\sum_{j=n+1}^{n_i} a_j\ge c_1 (n_i-n)
$
for every $ 0 \le n < n_i $ and  $i=1, \dots, l$. \fle \dem See \cite[Lemma 3.1]{ABV}. \cqd

We say that the {\em frequency of $(\sigma,\delta)$-hyperbolic times}\index{frequency} for $x\in M$ is bigger
than $\theta>0$ if, for large $n\in\NN$, there are $\ell \ge\th n$ and integers $1\le n_1<n_2\dots <n_\ell\le n$
which are $(\sigma,\delta)$-hyperbolic times for $x$.

\cpr \label{pr.hyperbolic1} Assume that $f\colon M\to M$ is
non-uniformly expanding. Then there are $0<\sigma<1$, $\delta>0$
and $\theta>0$ (depending only on $\lambda$ and on the derivative
of $f$) such that the frequency of $(\sigma,\delta)$-hyperbolic
times for Lebesgue almost all $x\in M$ is bigger than $\theta$.
 \fpr
  \dem
Assuming that \eqref{NUE} holds for $x\in M$, then for large $N\in\NN$ we have
\begin{equation*}
\sum_{j=0}^{N-1} -\log\|Df(f^j(x))^{-1}\| \ge \lambda N\,.
\end{equation*}
Take $\beta>0$ given by Definition~\ref{d.nd}, and fix any $\rho>\beta$. Then (s$_2$) implies that
\begin{equation}
\label{e.bound2} \left| \,\log \|Df(x)^{-1}\|\,\right| \le \rho \, |\log \dist(x,\cc)|
\end{equation}
for every $x$ in a neighborhood $V$ of $\cc$. Fix $\varepsilon_1>0$ so that $\rho\varepsilon_1 \le \lambda/2$,
and let $r_1>0$ be  so that
\begin{equation}
\label{e.bound1} \sum_{j=0}^{N-1}\log\dist_{r_1}(f^j(x),\cc)\ge -\varepsilon_1 N\,.
\end{equation}
The assumption of slow recurrence to the critical set ensures that this is possible. Fix any $K_1 \ge \rho \,
|\log r_1|$ large enough so that it is also an upper bound for $-\log\|Df^{-1}\|$ on the complement of $V$. Then
let $J$ be the subset of times $1 \le j \le N$ such that $-\log\|Df(f^{j-1}(x))^{-1}\| > K_1$, and define
$$
a_j=\left\{\begin{array}{ll}
           -\log\|Df(f^{j-1}(x))^{-1}\| & \text{if } j\notin J \\
           0 & \text{if } j\in J.
           \end{array}\right.
$$
By construction, $a_j\le K_1$ for $1\le j \le N$. Note that if $j\in J$ then $f^{j-1}(x)\in V$. Moreover, for
each $j\in J$
$$
\rho \, |\log r_1| \le K_1 < -\log\|Df(f^{j-1}(x))^{-1}\| < \rho \, |\log \dist(f^{j-1}(x),\cc)|,
$$
which shows that  $\dist(f^{j-1}(x),\cc) < r_1$ for every $j\in J$. In particular,
$$\dist_{r_1}(f^{j-1}(x),\cc)=\dist(f^{j-1}(x),\cc)<r_1,\quad\text{for all
$j\in J$}.$$ Therefore, by (\ref{e.bound2}) and (\ref{e.bound1}),
$$
\sum_{j\in J} -\log\|Df(f^{j-1}(x))^{-1}\| \le \rho \sum_{j\in J} | \log \dist(f^{j-1}(x),\cc) | \le \rho \,
\varepsilon_1 N.
$$
We have chosen $\varepsilon_1>0$ in such a way that the last term is less than $\lambda N/2$. As a consequence,
$$
\sum_{j=1}^{N} a_j =  \sum_{j=1}^{N} -\log\|Df(f^{j-1}(x))^{-1}\| -
                      \sum_{j\in J} -\log\|Df(f^{j-1}(x))^{-1}\|
                   \ge \frac{\lambda}{2} N \,.
$$
Thus, we have checked that we may apply Lemma~\ref{l.pliss} to the numbers $a_1,\dots,a_N$, with
$c_1=\lambda/4$, $c_2=\lambda /2$, and $A=K_1$. The lemma provides $\theta_1>0$ and $l_1\ge \theta_1 N$ times $1
\le p_1< \cdots <p_{l_1}\le N$ such that
\begin{equation}
\label{e.conclusion1} \sum_{j=n+1}^{p_i} -\log\|Df(f^{j-1}(x))^{-1}\| \ge \sum_{j=n+1}^{p_i} a_j \ge
\frac\lambda4(p_i-n) \,
\end{equation}
for every $0 \le n < p_i$ and $1\le i\le l_1$.


Now fix $\varepsilon_2>0$ small enough so that $\varepsilon_2 < \theta_1b \lambda/4$, and let $r_2>0$ be such
that
\begin{equation}\label{e.e2r2}
\sum_{j=0}^{N-1} \log\dist_{r_2}(f^j(x),\cc)\ge - \varepsilon_2 N \,.
\end{equation}
Let $c_1=- b \lambda/4$, $c_2=-\varepsilon_2$, $A=0$, and
$$
\theta_2=\frac{c_2-c_1}{A-c_1}=1-\frac{4\varepsilon_2}{b\, \lambda}\,.
$$
Applying Lemma~\ref{l.pliss} to $a_j=\log\dist_{r_2}(f^{j-1}(x),\cc)$, with $1\le j\le N$, we conclude that
there are $l_2\ge\theta_2 N$ times $1\le q_1 < \cdots < q_{l_2}\le N$ such that
\begin{equation}
\label{e.conclusion2} \sum_{j=n}^{q_i-1} \log\dist_{r_2}(f^j(x),\cc) \ge -\frac{b \lambda}4 \, (q_i - n)
\end{equation}
for every $0 \le n < q_i$ and $1\le i \le l_2$\,.

Finally, our condition on $\varepsilon_2$ means that $\theta_1+\theta_2>1$. Let $\theta=\theta_1+\theta_2-1$.
Then there exist $l=(l_1+l_2-N) \ge \theta N$ times $1 \le n_1< \cdots < n_l\le N$ at which
(\ref{e.conclusion1}) and (\ref{e.conclusion2}) occur simultaneously:
$$
\sum_{j=n}^{n_i-1} -\log\|Df(f^j(x))^{-1}\| \ge \frac\lambda4(n_i-n)
$$
and
$$
\sum_{j=n}^{n_i-1} \log\dist_{r_2}(f^j(x),\cc) \ge -\frac{b \lambda}4 (n_i - n),
$$
for every $0 \le n < n_i$ and $1\le i \le l$. Letting $\sigma=e^{-\lambda/4}$ we easily obtain from the
inequalities above
$$
\prod_{j=n_i-k}^{n_i-1} \|Df(f^j(x))^{-1}\| \le \sigma^{k} \quad\text{and}\quad \dist_{r_2}(f^{n_i-k}(x),\cc)
\ge \sigma^{b k}
$$
for every  $1\le i \le l$ and $1\le k\le n_i$. In other words, all those $n_i$ are $(\sigma,\delta)$-hyperbolic
times for $x$, with $\delta=r_2$. \cqd

\cre\label{r.strong} From the proof of the previous proposition one easily sees that condition
\eqref{e.faraway1} in the definition of non-uniformly expanding map is not needed in all its strength for the
proof work. Actually, we have only  used \eqref{e.faraway1} in \eqref{e.bound1} and \eqref{e.e2r2}. Hence, it is
enough that \eqref{e.faraway1} holds for $\vare=\min\{\vare_1,\vare_2\}$ and
 $\delta=\max\{r_1,r_2\}$. \fre

 \cre\label{r.start} Observe that the proof of
Proposition~\ref{pr.hyperbolic1} also gives that if for some $x\in M$ and $N\in \NN$
\begin{equation*}
\sum_{j=0}^{N-1} -\log\|Df(f^j(x))^{-1}\| \ge \lambda N \qand \sum_{j=0}^{N-1} \log\dist_{\delta}(f^j(x),\cc)\ge
- \varepsilon N \,
\end{equation*}
(where $\vare$ and $\delta$ chosen as in Remark~\ref{r.strong}), then there exist $1 \le n_1< \cdots < n_l\le N$
with $l\ge \theta N$ such that $n_i$ is a $(\sigma,\delta)$-hyperbolic time for $x$
for every $1\le i \le l$.  \fre

 \cco \label{co.distortion0} There exists $C_0=C_0(B,\beta,b,\sigma)>0$ such that for
every hyperbolic pre-ball $V_n$ and every $y, z\in V_n$
 $$ \log\frac{|\det Df^n (y)|}{|\det Df^n (z)|}
 \le C_0\dist(f^{n}(y),f^{n}(z)).
 $$
  \fco

\dem It suffices to take $C_0 \ge \sum_{k=1}^{\infty} 2^\beta B \sigma^{(1/2-b\beta)k}$; recall that $b\beta <
1/2$. \cqd

\cco \label{co.distortion}  There exists $C_1=C_1(C_0)>0$  such that for every hyperbolic pre-ball $V_n$ and
every $y, z\in V_n$
 $$  \frac{1}{C_1} \le
 \frac{|\det Df^n
(y)|}{|\det Df^n (z)|}\le C_1 \,.
 $$
 \fco

 \dem Take
 $C_1=\exp(C_0D),$ where $D$ is the diameter of $M$.
 \cqd

We finish this section deriving an useful consequence of the existence of positive frequency of hyperbolic
times.

\cle\label{c:hyperbolic3} Let $A\subset M$ be a set with positive Lebesgue measure whose points have frequency
of $(\sigma,\delta)$-hyperbolic times bigger than $\theta>0$. Then there is $n_0\in\NN$ such that for $n\ge n_0$
    $$
   \frac{1}{n}\sum_{j=1}^{n}\frac{\leb(A\cap H_{j})}{\leb(A)} \geq
   \frac\theta2,
   $$
  where $H_j$ is the
 set of points that have  $j$ as a
 $(\sigma,\delta)$-hyperbolic
 time.
\fle \dem Since we are assuming that points in $A$  have frequency of $(\sigma,\delta)$-hyperbolic times bigger
than $\theta>0$, then there are $n_0\in\NN$ and a  set $B\subset A$ with $m(B)\ge m(A)/2$ such that for every
$x\in B$ and
 $n\ge n_0$ there are $(\sigma,\delta)$-hyperbolic times
$0<n_1<n_2<\dots <n_\ell\le n$  for $x$ with $\ell \ge\th n$.
 Take now $n\ge n_0$ and let $\xi_n$ be the measure in $\{1,\dots,n\}$ defined
by $\xi_n(J)=\# J/n$, for each subset $J$. Then, using Fubini's Theorem
\begin{eqnarray*}
\frac{1}{n} \sum_{j=1}^{n}\leb(B\cap H_j)
& =& \int \left(\int_B \Chi(x,i)\,d\leb(x)\right)d\xi_n(i) \\
& = & \int_B \left(\int \Chi(x,i)\,d\xi_n(i)\right)d\leb(x),
\end{eqnarray*}
where $\Chi(x,i)=1$ if $x\in H_i$, and $\Chi(x,i)=0$ otherwise. Since for every $x\in B$ and
 $n\ge n_0$ there are $0<n_1<n_2<\dots <n_\ell\le n$ with $\ell \ge\th n$
 such that $x\in H_{n_i}$ for $1\le i\le\ell$, then the
integral with respect to $d\xi_n$ is larger than $\theta$. So, the last expression in the formula above is
bounded from below by $\theta \leb(B)\ge\theta \leb(A)/2$. \cqd


\section{Markov structures}\label{ch.rates}

 The aim of this section is to  show that non-uniformly expanding
transformations induce piecewise expanding Markovian return maps. This has been proved in \cite{ALP} and we
follow the proof therein. Detailed  proofs of most results are presented here in order to show how constants
depend on one another.



\cte\label{t:Markov towers}
 Let $ f:M\to M $ be a $C^2$
non-uniformly expanding transitive transformation. Then $f$
induces some piecewise expanding Markovian return map on a ball
$\Delta\subset M$. Moreover, if there exist $C, \gamma>0 $ such
that
    $
    \leb(\Gamma_n) \leq Cn^{-\gamma},
    $
 then  there is $C'>0$ such that the return time function satisfies
    $ \leb\{R>n\}\leq
    C'n^{-\gamma}.  $
   \fte

Assuming that $f$ is a non-uniformly expanding map, then by Proposition~\ref{pr.hyperbolic1} there are $\sigma$,
$\delta$ and $\theta$ such that Lebesgue almost every $x\in M$ has frequency of $(\si,\de)$-hyperbolic times
greater than $\theta$. From the transitivity of $f$ and by \cite[Lemma 2.5]{ALP} we may fix  $ p \in M $ and
$N_0\in\NN$ for which
\begin{equation}\label{defN0}
\cup_{j=0}^{N_0}f^{-j}\{p\}\,\text{ is $\delta_1/3$-dense in $ M $ and disjoint from $ \cc $,}
 \end{equation}
  where
$\delta_1>0$ is the radius of hyperbolic balls given by
Proposition~\ref{pr.distorcao}. Take constants $\vare>0$ and
$\delta_0>0$ so that
 $$
 \sqrt\delta_0 \ll \delta_1/2\qand 0< \varepsilon \ll \delta_0
 .
 $$
Let us introduce a couple of auxiliary lemmas.

\begin{Lemma}\label{le:grow2}
There are constants $ K_0,  D_{0} 
>0$ depending only on $ f$, $\sigma$, $\delta_1 $ and the point~$
p $, such that for any ball $ B\subset M $ of radius $ \delta_1 $
there are an open set $ V\subset B $ and an integer $ 0\leq m \leq
 N_{0} $ for which:
   \begin{enumerate}
   \item $ f^{m} $ maps $V$ diffeomorphically onto $B(p,2\sqrt\delta_{0}) $;
   \item
   for each $x,y\in V$
   $$
\log\left|\frac{\det Df^{m}(x)}{\det Df^{m}(y)}\right|
 \le
D_0 \dist(f^{m}(x),f^{m}(y)).
$$
\end{enumerate}
Moreover, for each  $ 0\leq j \leq
 N_{0} $  the $j$-preimages of $B(p,2\sqrt\delta_{0})$
       are all disjoint from $\cc$, and
  for $x$
       belonging to any such $j$-preimage we have
 $
 {K_0}^{-1}\le \|Df^j(x)\|
 \le K_0.
 $
\end{Lemma}

\begin{proof}
Since $\cup_{j=0}^{N_0}f^{-j}\{p\}$ is $ \delta_1/3 $ dense in $ M $ and disjoint from $\cc$, choosing $
\delta_{0}>0 $ sufficiently small we have that each connected component of the preimages of $B(p,
2\sqrt\delta_{0}) $ up to time $ N_{0} $ are bounded away from the critical set $ \cc $ and are contained in a
ball of radius $ \delta_1/3 $. This immediately implies that any ball $ B \subset M $ of radius $ \delta_1 $
contains a preimage $ V $ of $ B(p, 2\sqrt\delta_{0}) $ which is mapped diffeomorphically onto $ B(p,
2\sqrt\delta_{0}) $ in at most $ N_{0}$ iterates. Moreover, since the number of iterations and the distance to
the critical region are uniformly bounded, the volume distortion is uniformly bounded.

Observe that $ \delta_{0} $ and $N_0$ have been chosen in such a
way that all the connected components of the preimages of $B(p,
2\sqrt\delta_{0})$ up to time $ N_{0} $  are uniformly bounded
away from the critical set $ \cc $, and so there is some constant
$K_0>1$ 
such that
 $
 {K_0}^{-1}\le \|Df^m(x)\|
 \le K_0
 $
 for all $1\le m\le N_0$ and $x$ belonging to an $m$-preimage of $B(p,
2\sqrt\delta_{0})$.
\end{proof}


\cle\label{le:grow1}
    There exists $  N_{\varepsilon}
    > 0 $    such that any ball $ B
    \subset M $ of radius $ \varepsilon$ contains a hyperbolic pre-ball $
    V_{n}\subset B $ with $ n\leq  N_{\varepsilon} $.
\fle

\begin{proof}
    Take any $ \varepsilon > 0 $ and a ball $B(z,\vare)$. By
    Proposition~\ref{pr.distorcao} we may choose $n_{\varepsilon}\in\NN$
    large enough so that any hyperbolic pre-ball $ V_{n} $ associated to
    a
    hyperbolic time $ n\geq  n_{\varepsilon} $ has diameter not exceeding
     $\vare/2$.  Now
    notice that by Proposition~\ref{pr.hyperbolic1} Lebesgue almost every point has an infinite number of hyperbolic times and
    therefore
    $$
    \leb\left(M\setminus \mcup_{j= n_{\varepsilon}}^{n}H_{j}\right) \to 0
    \quad\text{ as } n\to\infty.
    $$
    Hence, it is possible to choose $N_{\varepsilon}\in \NN$ such
    that
    $$
    \leb\left(M\setminus \mcup_{j=
    n_{\varepsilon}}^{N_\vare}H_{j}\right)
   < m(B(z, \varepsilon/2)).
    $$
  This ensures that there is
    a point $ \hat x \in B(z, \varepsilon/2) $ with a hyperbolic time
    $ n\leq  N_{\varepsilon} $ and associated hyperbolic pre-ball $
    V_{n}(x)$ contained in $ B(z, \varepsilon) $.
\end{proof}

\cre\label{re.Ne} Observe that if $n$ is a hyperbolic time for $f$, then $n$ is also a hyperbolic time for every
map in a sufficiently small $C^1$ neighborhood of $f$. Hence, for given $\vare>0$ the integer $N_\vare$ may be
taken uniform in a whole $C^1$ neighborhood of $f$, and  only depending on $\vare$, $\sigma$ and $\delta_1$.

\fre


\subsection{The partitioning algorithm}\label{s.algoritmo}

Here we describe the construction of the  partition (mod 0) of $
\Delta_{0}=B(p,\delta_{0}) $.  
We introduce
 neighborhoods of $ p $
$$\Delta^{0}_{0}=
B(p,\delta_{0}),\quad \Delta^{1}_{0}=B(p,2\delta_{0}),\quad
\Delta^{2}_{0}=B(p,\sqrt\delta_{0})\qand
\Delta^{3}_{0}=B(p,2\sqrt\delta_{0}).$$ For $0<\sigma<1$ given by
Proposition~\ref{pr.hyperbolic1}, let
$$
I_{k}=\left\{x\in\Delta^{1}_{0}\: : \:\delta_{0}(1+\sigma^{k/2}) <
\dist(x,p) < \delta_{0}(1+\sigma^{(k-1)/2})\right\},\quad  k\ge 1,
$$
be a partition (mod 0) into countably many rings of $ \Delta_0^{1}\setminus \Delta_0 $. The construction of the
partition of $\Delta_0$ is inductive and we describe precisely the general step of the induction below. 

 Take $ R_{0} $ 
some large integer to be determined latter; we ignore any dynamics
occurring up to time $ R_{0} $. 
%
%
%
Assume that sets $\Delta_{i}$, $A_{i}$, $A_{i}^\vare$
 $B_{i}$, $\{R=i\}$ and functions $ t_{i}:
\Delta_{i}\to\mathbb N $ are defined for all $ i\leq n-1 $. For $
i\leq R_{0} $ we just let $ A_{i}=A_{i}^\vare=\Delta_{i}=
\Delta_{0}$, $ B_{i}=\{R=i\}=\emptyset $ and $ t_{i}\equiv 0 $.
Now let $(U_{n,j}^3)_j$ be the connected components of
 $
 f^{-n}(\Delta_0)\cap A_{n-1}^\vare
 $
 contained in hyperbolic pre-balls $V_m$, with $n-  N_0\le
 m\le n$, which are mapped onto $\Delta_0^3$ by $f^n$.
 Take
$$
U_{n,j}^i=U_{n,j}^3\cap f^{-n}\Delta_0^i,\quad i=0,1,2,
$$
and set $R(x)=n$ for $x\in U_{n,j}^0$. Take also
$$
\Delta_{n}=\Delta_{n-1}\setminus \{R=n\}.
$$
The definition of the function $ t_{n}:\Delta_{n}\to \mathbb N $
is slightly different in the general case:
\begin{equation*}
    t_{n}(x) =
    \begin{cases}
    s & \text{ if } x\in U_{n,j}^{1}\setminus U_{n,j}^{0} \text{ and }
    f^{n}(x) \in I_{s} \text{ for some $j$,} \\
    0 & \text{ if } x\in A_{n-1} \setminus \bigcup_{j} U^{1}_{n,j},\\
    t_{n-1}(x)-1 & \text{ if } x\in B_{n-1}\setminus \bigcup_{j} U^{1}_{n,j}.
    \end{cases}
\end{equation*}
Finally let
$$
A_{n}= \{x\in\Delta_{n}: t_{n}(x) = 0\}, \quad B_{n}=
\{x\in\Delta_{n}: t_{n}(x) > 0\}
$$
and
 $$
 A_{n}^\vare= \{x\in\Delta_{n}:
 \dist(f^{n+1}(x),f^{n+1}(A_n))<\vare\}.
 $$
At this point we have completely described the inductive
construction of the sets  $A_n$, $A_n^\vare$,  $B_n$ and
$\{R=n\}$.

\medskip

 The construction detailed before provides an
 algorithm for the definition of a family of topological balls
 contained in $ \Delta_0 $ and satisfying the Markov property as required.
This algorithm does indeed
 produce a partition mod 0 of $ \Delta_0 $; see \cite[Lemma 3.1]{ALP}.

Associated to each component $U^0_{n-k}$ of $\{R=n-k\}$, for some $k>0$, we have a collar $U^1_{n-k}\setminus
U^0_{n-k}$ around it; knowing that the new components of $\{R=n\}$ do not intersect ``too much"
$U^1_{n-k}\setminus U^0_{n-k}$ is important for preventing overlaps on sets of the partition.  This is indeed
the case as long as $\vare>0$ is taken small enough.

\cle\label{l.claim} If $\vare>0$ is sufficiently small, then $U_{n}^1\cap\{t_{n-1}\ge 1\}=\emptyset$ for each
$U^1_n$. \fle

\dem Take some $k>0$ and let  $U_{n-k}^0$ be a component of $\{R=n-k\}$. Let $Q_{k}$ be the part of $U_{n-k}^1$
that is mapped by $f^{n-k}$ onto $I_{k}$ and assume that $Q_{k}$ intersects some $U_{n}^3$. Recall that, by
construction, $Q_{k}$ is precisely the part of $U_{n-k}^1$ on which $t_{n-1}$ takes the value 1. Letting $q_1$
and $q_2$ be any two points in distinct components (inner and outer) of the boundary of $Q_{k}$, we have
 by Proposition~\ref{pr.distorcao}  and Lemma~\ref{le:grow2}
 \begin{equation}\label{e.zq1}
 \dist(f^{n-k}(q_1),f^{n-k}(q_2))\le
 K_0\sigma^{(k-N_0)/2}\dist(f^{n}(q_1),f^{n}(q_2)).
 \end{equation}
We also have
 \begin{eqnarray*}
\dist(f^{n-k}(q_1),f^{n-k}(q_2))&\ge& \delta_{0}(1+\sigma^{(k-1)/2})-\delta_{0}(1+\sigma^{k/2})\\&=& \delta_{0}
\sigma^{k/2}(\sigma^{-1/2}-1),
 \end{eqnarray*}
which combined with (\ref{e.zq1})  gives
$$
 \dist(f^{n}(q_1),f^{n}(q_2))\ge K_0^{-1}\sigma^{N_0/2}\delta_{0}(\sigma^{-1/2}-1).
 $$
On the other hand, since $ U^{3}_{n}\subset A_{n-1}^{\varepsilon}$ by construction of $U^{3}_{n}$, taking
 \begin{equation}\label{eq.vare}
 \vare<K_0^{-1}\sigma^{N_0/2}\delta_{0} (\sigma^{-1/2}-1)
 \end{equation}
we  have $U_{n}^3\cap\{t_{n-1}>1\}=\emptyset$. This implies $U_{n}^1\cap\{t_{n-1}\ge 1\}=\emptyset$. \cqd


  \subsection{Expansion. }\label{sub.exp}

Recall that by construction, the return time $ R $ for an element $ U $ of the partition $ \mathcal P $ of $
\Delta_0 $ is formed by a certain number $ n $ of iterations given by the hyperbolic time of a hyperbolic
pre-ball $ V_{n}\supset U $, and a certain number $ m\leq N_{0} $ of additional iterates which is the time it
takes to go from $ f^{n}(V_{n}) $ which could be anywhere in $ M $, to $ f^{n+m}(V_{n}) $ which covers $
\Delta_0 $ completely. It follows from Proposition~\ref{pr.distorcao} and Lemma~\ref{le:grow2} that
$$
\|Df^{n+m}(x)^{-1}\|\le
\|Df^{m}(f^n(x))^{-1}\|\cdot\|Df^{n}(x)^{-1}\|<K_0\sigma^{n/2}\le
K_0\sigma^{(R_0-N_0)/2}.
$$
By taking $ R_{0} $ sufficiently large we can make this last
expression smaller than 1.


\subsection{Bounded distortion.}\label{sub.bd}
For the  bounded distortion estimate  in Definition~\ref{de.expanding} we need to show that there exists a
constant $ K> 0 $ such that for any $ x, y $ belonging to an element $ U\in\mathcal P $ with return time $ R $,
we have
$$
\log\left|\frac{\det Df^{R}(x)}{\det Df^{R}(y)}\right| \leq K
\dist(f^{R}(x), f^{R}(y)).
$$
Recall that by construction, the return time $ R$ for an element $
U $ of the partition $ \mathcal P $ of $ \Delta_0 $ is formed by a
certain number $ n $ of iterations given by the hyperbolic time of
a hyperbolic pre-ball $ V_{n}\supset U $, and a certain number $
m=R-n\leq N_{0} $ of additional iterates which is the time it
takes to go from $ f^{n}(V_{n}) $ to $ \Delta_0 $  and cover it
completely.
 By the chain rule
    \begin{align*}
\log\left|\frac{\det Df^{R}(x)}{\det Df^{R}(y)}\right| & =
\log\left|\frac{\det Df^{R-n}(f^n(x))}{\det
Df^{R-n}(f^n(y))}\right| + \log\left|\frac{\det Df^{n}(x)}{\det
Df^{n}(y)}\right|.
    \end{align*}
For the first term in this last sum we observe that by
Lemma~\ref{le:grow2} we have
$$
\log\left|\frac{\det Df^{R-n}(f^n(x))}{\det
Df^{R-n}(f^n(y))}\right|
 \le
D_0 \dist(f^{R}(x),f^{R}(y)).
$$
For the second term in the sum above, we may apply
Corollary~\ref{co.distortion0} and obtain
 $$
\log\left|\frac{\det Df^{n}(x)}{\det Df^{n}(y)}\right| \le
 C_0\dist(f^{n}(x),f^{n}(y))
 .$$
Also by Lemma~\ref{le:grow2} we may write
 $$\dist(f^{n}(x),f^{n}(y))\le K_0 \dist(f^{R}(x),f^{R}(y)).$$
Thus we just have to take $K=D_0+C_0K_0$.

In the next lemma we show that the bounded distortion condition in \cite{AV} is satisfied in our context.

\cle\label{le.bd} For each $U\in \cp$ we have
$$
 \frac{\left\|
D\left(J\circ (F\vert_U)^{-1}\right)\right\|} {\left|J\circ
(F\vert_U)^{-1}\right|}<K,
 $$
 where $J=\det DF$ is the Jacobian of $F$.
 \fle
 \dem For simplicity we assume $\Delta\subset\RR^d$.
Observe that
 $$
 \frac{\left\|D\left(J\circ (F\vert_U)^{-1}\right)\right\|}
              {\left|(J\circ (F\vert_U)^{-1})\right|}
 =\left\|D\left(\log\left| J\circ (F\vert_U)^{-1}\right|\right) \right\|.
 $$
 Thus we just have to prove that the   functions $\log\left| J\circ
 (F\vert_U)^{-1}\right|$, $U\in\cp$,  have  derivatives uniformly bounded by $K$.
Take any point $x$ in the interior of $ \Delta$ and $v$ a vector
of the canonical basis of $\RR^d$. By the bounded distortion
condition of Definition~\ref{de.expanding} we have for small
$t\in\RR$
 \begin{align*}
 \log\left| J\circ (F\vert_U)^{-1}\right|(x+tv)&
 -\log\left| J\circ (F\vert_U)^{-1}\right|(x)\\
 &\le K\dist (F((F\vert_U)^{-1}(x+tv)),F((F\vert_U)^{-1}(x)))\\&= Kt.
 \end{align*}
 This implies the uniform bound on derivatives that we need.
 \cqd

 \subsection{Metric estimates}\label{s.check}
Now we prove that the construction performed above does indeed
produce a partition of $\Delta_0$ as in the Theorem~\ref{t:Markov
towers}, modulo a zero Lebesgue measure subset.  We split our
argument into two parts.

\medskip

\subsubsection{Estimates derived from the construction} In this first part we obtain some estimates
relating the Lebesgue measure of the sets $A_n$, $B_n$ and
$\{R>n\}$ with the help of specific information extracted from the
inductive construction we performed in
Subsection~\ref{s.algoritmo}.

 \cle\label{l.flowb} There exists a constant $a_0>0$ (not depending on $\delta_0$)
  such that  $$\leb(B_{n-1}\cap A_n)\ge a_0\leb(B_{n-1})$$ for every
$n\ge1$.
 \fle
  \dem
  It is enough to see that this holds for each connected component of $B_{n-1}$
 at a time. Let  $C$ be a component of $B_{n-1}$ and $Q$
be  its outer ring corresponding to $t_{n-1}=1$. Observe that by Lemma~\ref{l.claim} we have $Q=C\cap A_n$.
Moreover, there must be some $k<n$ and a component $U^0_k$ of $\{R=k\}$ such that $f^k$ maps $C$
diffeomorphically onto $\bigcup_{i=k}^\infty I_i$ and $Q$ onto $I_k$, both with distortion bounded by $C_1$ and
$e^{D_0L}$,
where $L$ is the diameter of $M$; 
 cf. Corollary~\ref{co.distortion}
and Lemma~\ref{le:grow2}.
 Thus, it is sufficient to compare the Lebesgue measures of $\bigcup_{i=k}^\infty
I_i$ and $ I_k$. We have
 $$
 \frac{\leb( I_k)}{\leb(\bigcup_{i=k}^\infty
I_i)}\thickapprox\frac{
[\delta_0(1+\sigma^{(k-1)/2})]^d-[\delta_0(1+\sigma^{k/2})]^d}
{[\delta_0(1+\sigma^{(k-1)/2})]^d-\delta_0^d}\thickapprox
1-\sigma^{1/2}.
 $$
 Clearly this proportion does not depend on  $\delta_0$.
  \cqd

  \cle\label{l.flowa}  There exist $b_0,c_0>0$  with $b_0+c_0<1$
  such that for every $n\ge1$
\begin{enumerate}
 \item $\leb(A_{n-1}\cap B_n)\le b_0\leb(A_{n-1})$;
 \item $\leb(A_{n-1}\cap \{R=n\})\le c_0\leb(A_{n-1})$.
\end{enumerate}
Moreover $b_0\to 0$ and $c_0\to 0$  as $\delta_0\to 0$.
 \fle
\dem It is enough to prove these estimates for each neighborhood
of a component $U^0_n$ of $\{R=n\}$. Observe that by construction
we have $U^3_n\subset A_{n-1}^\vare$, which means that $U^2_n
\subset A_{n-1}$, because $\vare<\delta_0<\sqrt\delta_0$. Using
the  distortion bounds of $f^n$ on $U^3_n$ given by
Corollary~\ref{co.distortion} and Lemma~\ref{le:grow2}  we obtain
 $$
 \frac{\leb(U^1_n\setminus U^0_n)}{\leb(U^2_n\setminus U^1_n)}
 \thickapprox
 \frac{\leb(\Delta^1_0\setminus \Delta^0_0)}{\leb(\Delta^2_0\setminus \Delta^1_0)}
 \thickapprox
 \frac{\delta_0^d}{\delta_0^{d/2}}\ll 1,
 $$
which gives the first estimate. Moreover,
$$
 \frac{\leb( U^0_n)}{\leb(U^2_n\setminus U^1_n)}
 \thickapprox
 \frac{\leb( \Delta^0_0)}{\leb(\Delta^2_0\setminus \Delta^1_0)}
 \thickapprox
 \frac{\delta_0^d}{\delta_0^{d/2}}\ll 1,
 $$
and this gives the second one.
 \cqd

The next result asserts that a fixed proportion of $A_{n-1}\cap
H_n$ gives rise to new elements of the partition within a finite
number of steps (not depending on $n$).

    \cpr\label{p.construction}
    There exist $c_1>0$ and a positive integer $N=N(\vare)$ such that
    $$
     \leb\left(\mcup_{i=0}^N\big\{R=n+i\big\}\right)\ge c_1 \leb(A_{n-1}\cap H_{n})
    $$ for every $ n\ge1$.
    \fpr
\dem Take $r=5\delta_0K_0^{N_0}$, where $N_0$ and $K_0$ are given
by Lemma~\ref{le:grow2}. Let $\{z_j\}$ be a maximal set in
$f^n(A_{n-1}\cap H_n)$ with the property that $B(z_j,r)$ are
pairwise disjoint. By maximality we have
 $
 \mcup_j B(z_j,2r)\supset f^n(A_{n-1} \cap H_n).
 $
Let $x_j$ be a point in $H_n$ such that $f^n(x_j)=z_j$ and
consider the hyperbolic pre-ball $V_n(x_j)$ associated to $x_j$.
Observe that $f^n$ sends $ V_n(x_j)$ diffeomorphically onto a ball
of radius $\delta_1$ around $z_j$ as in
Proposition~\ref{pr.hyperbolic1}. In what follows, given $B\subset
B(z_j,\delta_1)$, we will simply denote
 $(f^{n}\vert V_n(x_j))^{-1}(B)$ by $f^{-n}(B)$.

Our aim now is to prove that $f^{-n}(B(z_j, r))$ contains some
component of $\{R=n+k_j\}$  with  $0\le k_j\le N_\vare+N_0$. We
start by showing that
 \begin{equation}\label{eq.tnk}
  t_{n+k_j}\vert f^{-n}(B(z_j,\vare))>0\quad\text{for some $0\le k_j\le
  N_\vare+N_0$}.
  \end{equation}
Assume by contradiction that
  $
  t_{n+k_j}\vert f^{-n}(B(z_j,\vare))=0$ for all $0\le k_j\le
  N_\vare+N_0
  $.
 This implies that $f^{-n}(B(z_j,\vare))\subset A_{n+k_j}^\vare$
 for all $0\le k_j\le N_\vare+N_0$.
 Using Lemma~\ref{le:grow1} we may find  a hyperbolic pre-ball
 $ V_{m}\subset B(z_j,\vare) $ with $ m\leq  N_{\varepsilon} $.
 Now, since $f^m(V_{m})$ is a ball $B$ of radius $\delta_1$ it follows
from Lemma~\ref{le:grow2} that there is some $V\subset B$
and $m'\le N_0$ with $f^{m'}(V)=\Delta_0$. Thus, taking $k_j=m+m'$
we have that $0\le k_j\le N_\vare+N_0$ and $f^{-n}(V_m)$ is an
element of $\{R=n+k_j\}$ inside $f^{-n}(B(z_j,\vare))$. This
contradicts the fact that $t_{n+k_j}\vert f^{-n}(B(z_j,\vare))=0$
 for all $0\le k_j\le N_\vare+N_0$, and so (\ref{eq.tnk}) holds.

Let $k_j$ be the
 smallest integer $0\le k_j\le N_\vare+N_0$ for which
$t_{n+k_j}\vert f^{-n}(B(z_j,\vare))>0$.
Since
  $
  f^{-n}(B(z_j,\vare))\subset A_{n-1}^\vare\subset \{t_{n-1}\le1 \},
  $
 there must be some element
$U^{0}_{n+k_j}(j)$ of $\{R=n+k_j\}$ for which
 $
 f^{-n}(B(z_j,\vare))\cap U_{n+k_j}^1(j)\neq\emptyset.
 $
Recall that by definition $f^{n+k_j}$ sends $U_{n+k_j}^1(j)$ diffeomorphically onto $\Delta_0^1$, the ball  of
radius $(1+s)\delta_0$ around $p$. From time $n$ to $n+k_j$ we may have some final ``bad" period of length at
most $N_0$ where the derivative of $f$ may contract, however being bounded from below by $1/K_0$ in each step.
Thus, the diameter of $f^n(U_{n+k_j}^1(j))$ is at most $4\delta_0K_0^{N_0}$. Since $B(z_j,\vare)$ intersects
$f^n(U_{n+k_j}^1(j))$ and $\vare<\delta_0<\delta_0K_0^{N_0}$, we have by the definition of $r$ that
 $
 f^{-n}(B(z_j,r))\supset U_{n+k_j}^0(j).
 $
Thus we have shown  that $f^{-n}(B(z_j, r))$ contains some
component of $\{R=n+k_j\}$ with  $0\le k_j\le N_\vare+N_0$.
Moreover, since $n$ is a hyperbolic time for $x_j$, we have by the
distortion control given by Corollary~\ref{co.distortion}
 \begin{equation}\label{eq.quo1}
 \frac{ \leb(f^{-n}(B(z_j,2r)))}{\leb(f^{-n}(B(z_j,r)))}
 \le
 {C_1}\frac{\leb(B(z_j,2r))}{\leb(B(z_j,r))}
 \end{equation}
and
 \begin{equation}
\label{eq.quo2} \frac{\leb(f^{-n}(B(z_j,r)))}{\leb(
U_{n+k_j}^0(j))}
 \le
 {C_0}\frac{\leb(B(z_j,r))}{\leb(f^n(U_{n+k_j}^0(j)))}.
 \end{equation}
Here we are implicitly assuming that
\begin{equation}
 \label{eq.r}
 r=r(\delta_0)<{\delta_1}/2.
\end{equation}
 This
can be done by taking $\delta_0$ small enough. Note that estimates
on $N_0$ and $K_0$ improve when we diminish $\delta_0$.

From time $n$ to time $n+k_j$ we have at most $k_j=m_1+m_2$
iterates with $m_1\le N_\vare$, $m_2\le N_0$ and
$f^n(U_{n+k_j}^0(j)))$ containing some point $w_j\in H_{m_1}$. By
the definition of $(\sigma,\delta)$-hyperbolic time we have
$\dist_\delta (f^i(x),\cc)\ge \sigma^{bN_\vare}$ for every $0\le
i\le m_1$, which 
implies that there is some constant $D=D(\vare)>0$ such that
$|\det (Df^i(x))|\le D$ for $0\le i\le m_1$ and $x\in
f^n(U_{n+k_j}^0(j))$. On the other hand, since the first $N_0$
preimages of $\Delta_0$ are uniformly bounded away from $\cc$ we
also have some $D'>0$ such that $|\det (Df^i(x))|\le D'$ for every
$0\le i\le m_2$ and $x$ belonging to an $i$ preimage of
$\Delta_0$. Hence,
$$\leb(f^n(U_{n+k_j}^0(j)))\ge \frac{1}{DD'}\leb(\Delta_0),$$
which combined with (\ref{eq.quo2}) gives
 $$
 \leb(f^{-n}(B(z_j,r)))\le C\leb(
U_{n+k_j}^0(j)),
 $$
with $C$ only depending on $C_1$, $D$, $D'$, $\delta_0$ and the
dimension of $M$. We also deduce from (\ref{eq.quo1}) that
 $$
 \leb(f^{-n}(B(z_j,2r)))\le C'\leb(
f^{-n}(B(z_j,r)))
 $$
with $C'$ only depending on $C_1$ and the dimension of $M$.
 Finally let us compare the Lebesgue measure of
the sets $\bigcup_{i=0}^N\big\{R=n+i\big\}$ and $A_{n-1}\cap
H_{n}$. We have
 $$
 \leb\big(A_{n-1}\cap
 H_n \big)\le \sum_j \leb(f^{-n}(B(z_j,2r)))
 \le C'
 \sum_j \leb(f^{-n}(B(z_j,r))).
  $$
On the other hand, by the disjointness of the balls $B(z_j,r)$ we
have
$$
 \sum_j \leb(f^{-n}(B(z_j,r)))\le C
 \sum_j \leb( U_{n+k_j}^0(j)) \le C
 \leb\left(\bigcup_{i=0}^N\big\{R=n+i\big\}\right).
  $$
We just have to take $c_1=(CC')^{-1}$. \cqd

\cre \label{re.c1}It follows from the choice of the constants $D$
and $D'$ (and so also $C$ and $C'$) that the constant $c_1$ only
depends on the constants $\sigma$, $b$, $N_\vare$, $C_1$ and
$N_0$. \fre

\subsubsection{General estimates}\label{se.tail}
For the time being  we have taken a disk $\Delta_0$ of radius
   $\delta_0>0$ around a point $p$
   and defined inductively the subsets
  $A_n$, $B_n$, $\{R=n\}$ and $\Delta_n$ which are related in
   the following way:
    $$\Delta_n=\Delta_0\setminus\{R\le n\}=A_n\dot\cup B_n.$$
Since we are dealing with a non-uniformly expanding map, we also
have defined for each $n\in \NN$ the set  $H_n$ of points that
have $n$ as a $(\si,\delta)$-hyperbolic time, and the tail of
expansion $\Gamma_n$  as in \eqref{tailset}. From the definition
of $\Gamma_n$, Remark~\ref{r.start} and Lemma~\ref{c:hyperbolic3}
we deduce:

\begin{enumerate}
\item[(m$_1$) ] there is $\theta>0$ such that for every $n\ge1$
and every $A\subset M\setminus\Gamma_n$ with $m(A)>0$
    \[
   \frac{1}{n}\sum_{j=1}^{n}\frac{\leb(A\cap H_{j})}{\leb(A)} \geq \theta.
   \]
\end{enumerate}
Moreover, we have proved in Lemma~\ref{l.flowb},
Lemma~\ref{l.flowa} and Proposition~\ref{p.construction} that the
following metric relations also hold:
\begin{enumerate}
\item[(m$_2$) ] there is $a_0>0$ (bounded away from 0 with $\delta_0$) such that for $n\ge1$
 $${\leb(B_{n-1}\cap A_n)}\ge a_0{\leb(B_{n-1})};$$
\item[(m$_3$) ]  there are $b_0,c_0>0$  with $b_0+c_0<1$ and
$b_0,c_0\to 0$
 as $\delta_0\to 0$, such
that for $n\ge1$
 $$\frac{\leb(A_{n-1}\cap B_n)}{\leb(A_{n-1})}\le b_0\qand
 \frac{\leb(A_{n-1}\cap \{R=n\})}{\leb(A_{n-1})}\le c_0;$$
\item[(m$_4$) ] there is $c_1>0$ and an integer $N\ge 0$ such that for $ n\ge1$
  $$ \leb\left(\mcup_{i=0}^N\big\{R=n+i\big\}\right)\ge c_1 \leb(A_{n-1}\cap
    H_{n}).$$
\end{enumerate}
In the inductive process of construction of the sets $A_n$, $B_n$,
$\{R=n\}$ and $\Delta_n$ we have fixed some large integer $R_0$,
being this the first step at which the construction began. Recall
that $A_n=\Delta_n=\Delta_0$ and $B_n=\{R=n\}=\emptyset$ for $n\le
R_0$. For technical reasons we will assume that
 \begin{equation}\label{e.constants}
 R_0>\max\left\{2(N+1),\frac{12}\theta\right\}.
 \end{equation}
Note that since $N$ and $\theta$ do not depend on $R_0$ this is
always possible.

\smallskip

This is the abstract setting under which we will be completing the
proof of Theorem~\ref{t:Markov towers}. From now on we will only
make use of the metric relations (m$_1$)-(m$_4$) and will not be
concerned with any other properties about these sets.

\cle\label{p.ab}
    There is  $a_1>0$, with $a_1\to 0$ as $\delta_0\to 0$, such that for all $n\ge1$
    $$
    \leb(B_{n}) \leq a_1\leb(A_{n}) .
    $$
\fle \dem   Let us just mention how the constant $a_1>0$ appears.
By (m$_3$)
 \begin{equation}\label{eq.eta}
 \leb(A_n\cap A_{n-1})\ge
 \eta\leb(A_{n-1}),
 \end{equation}
where $\eta=1-b_0-c_0$. Then we take
 \begin{equation}
 \label{e.waa1}
 \widehat a=\frac {b_0+c_0}{a_0}\qand
 a_1=\frac{(1+a_0)b_0+c_0}{a_0\eta}.
 \end{equation}
 The proof now follows exactly as in \cite[Proposition 5.4]{ALP}.
 \cqd


\cco\label{c.omega}
    There exists  $c_2>0$ such that for every $n\ge1$
    $$
    \leb(\Delta_{n}) \le c_2\leb(\Delta_{n+1}) .
    $$
\fco
 \dem
Using (m$_3$) we obtain
 $$
 \leb(\Delta_{n+1})\ge \leb(A_{n+1}) \ge (1-b_0-c_0)\leb(A_n).
 $$
On the other hand, by Lemma~\ref{p.ab},
 $$\leb(\Delta_n)=\leb(A_n)+\leb(B_n)\le
 (1+a_1^{-1})\leb(A_n).
 $$
It is enough to take $c_2=(1+a_1^{-1})/(1-b_0-c_0)$.
  \cqd

At this point we are able to definitely specify  the choice of
$\delta_0$. First of all, let us recall that the number $\theta$
in (m$_1$) does not depend on $\delta_0$. Assume that
 $\leb(\Gamma_n) \leq  Cn^{-\gamma},$ for some $C,\gamma>0$,
and pick $\alpha>0$ such that
 \begin{equation}\label{e.alpha}
  \alpha<\left(\frac\theta{12}\right)^{\gamma+1}.
  \end{equation}
Then we choose $\delta_0>0$ small enough so that
  \begin{equation}\label{e.a1}a_1<2\alpha.
  \end{equation}
This is possible because
     $a_1\to 0$ as $\delta_0\to 0$ by Lemma~\ref{p.ab}.

%
%
%
%
Since  $m(\Delta_n)=m(A_n)+m( B_n)$, we easily deduce from (m$_4$)
and Lemma~\ref{p.ab} that if we take
\begin{equation}\label{b1}
    b_1=\frac{c_1}{1+a_1},
\end{equation}
then
 $$
 \leb\left(\cup_{i=0}^N\big\{R=n+i\big\}\right)\ge b_1
 \frac{\leb(A_{n-1}\cap H_{n})}{\leb(A_{n-1})}\leb(\Delta_{n-1}).
 $$
This immediately implies that
 \begin{equation}\label{eq.recorre}
 \leb\left(\Delta_{n+N}\right)\le \left(1-b_1
 \frac{\leb(A_{n-1}\cap H_{n})}{\leb(A_{n-1})}\right)\leb(\Delta_{n-1}).
 \end{equation}
At this point we  obtained some recurrence relation for the
Lebesgue measure of the sets $\Delta_n$. 
Since $(\Delta_n)_n$ forms a decreasing sequence of sets we
finally have
 \begin{equation}\label{e.media}
 \leb\left(\Delta_{n+N}\right)\le \exp\left(-\frac{b_1}{N+1}\sum_{j=R_0}^{n}
 \frac{\leb(A_{j-1}\cap
 H_{j})}{\leb(A_{j-1})}\right)\leb(\Delta_{0}).
 \end{equation}
We will complete the proof of Theorem~\ref{t:Markov towers} by considering several different cases, according to
the behavior of the proportions ${\leb(A_{j-1}\cap H_{j})}/{\leb(A_{j-1})}$.
We  define for each $n\ge 1$
 $$
 E_n=\left\{ j\le n\colon \frac{\leb(A_{j-1}\cap
 H_j)}{\leb(A_{j-1})}<\alpha\right\},
 $$
and
 $$F=\left\{n\in\NN\colon \frac{\#
 E_n}{n}>1-\frac{\theta}{12}\right\}.$$

\cpr\label{p.final} Take any  $n\in F$ with $n\ge R_0$. If
$\leb(A_n)\ge 2 \leb(\Gamma_n)$, then there is some $0<k=k(n)<n$
for which
 ${\leb(A_n)}<\left( k/n\right)^\gamma{\leb(A_k)}.
 $
 \fpr

 \dem See \cite[Proposition 6.1]{ALP}. \cqd

Let us now complete the proof of  Theorem~\ref{t:Markov towers}. 
From Lemma~\ref{p.ab} we get
\begin{equation}\label{e.adelta}
\leb(\Delta_n)\le (1+a_1) \leb(A_n).
 \end{equation}
 Hence, it is
enough to derive the tail estimate of Theorem~\ref{t:Markov towers} for $\leb(A_n)$ in the place of
$\leb\{R>n\}=\leb(\Delta_n)$. Given any large integer $n$, we consider  the following two cases:
\begin{enumerate}
\item If $n\in\NN\setminus F$, then by (\ref{e.media}) and
Corollary~\ref{c.omega} we have
 \begin{equation}\label{e.f}\nonumber
 \leb(\Delta_n)\le c_2^N\exp\left({-\frac{b_1\theta\alpha}{12(N+1)}
 (n-R_0)}\right){\leb(\Delta_0)}.
 \end{equation}
\item If $n\in F$, then we distinguish the next two subcases:
\begin{enumerate}
\item If $\leb(A_n)< 2 \leb(\Gamma_n)$, then nothing has to be
done. \item If $\leb(A_n)\ge 2 \leb(\Gamma_n)$, then we apply
Proposition~\ref{p.final} and get some $k_1<n$ for which
$${\leb(A_n)}<\left(\frac {k_1}n\right)^\gamma {\leb(A_{k_1})}.
 $$
\end{enumerate}
\end{enumerate}
The only case we are left to consider is 2(b). In such case,
either $k_1$ is in situation~1 or 2(a), or by
Proposition~\ref{p.final} we can find $k_2<k_1$ for which
 $$
 {\leb(A_{k_1})}<\left(\frac {k_2}{k_1}\right)^\gamma {\leb(A_{k_2})}.
 $$
Arguing inductively we are able to show that there is a sequence
of integers $0<k_s<\cdots<k_1<n$ for which   one of the following
situations eventually  holds:
 \begin{enumerate}
\item[(A)]\quad $\displaystyle{\leb(A_{n})<\left(\frac
{k_s}{n}\right)^\gamma
c_2^N\exp\left({-\frac{b_1\theta\alpha}{12(N+1)}
 (k_s-R_0)}\right){\leb(\Delta_0)}}$.
\item[(B)]\quad $\displaystyle{\leb(A_{n})<\left(\frac
{k_s}{n}\right)^\gamma \leb(\Gamma_{k_s})}$. \item[(C)]\quad
$\displaystyle{\leb(A_{n})<\left(\frac {R_0}{n}\right)^\gamma
\leb(\Delta_0)}$.
 \end{enumerate}
 In all these three situations we arrive at the desired conclusion  of
Theorem~\ref{t:Markov towers}. Situation (C) corresponds to
falling in case 2(b) above successively until $k_s\le R_0$.


\section{Uniformness}
\label{s.uniformity}
%

Let us remark that the ball on which the piecewise expanding Markovian return map is defined may be taken the
same for every map belonging to a sufficiently small $C^2$ neighborhood of a map $f$ in a uniform family. In
fact, we have taken the ball $\Delta_0$ centered at a point $ p \in M $ which has been chosen in  \eqref{defN0}.
Since $\delta_1$ may be chosen the same for every $f$ in a uniform family, and the radius $\delta_0$ of the ball
$\Delta_0$ may be taken uniform in a neighborhood of $f$ (see Remark~\ref{re.uniform}), then the point $p$ and
$N_0$, and hence the ball $\Delta_0$, may be taken the same for every map belonging to a sufficiently small
$C^2$ neighborhood of $f$. Observe also that by an implicit function argument  the critical set varies
continuously with the map in the $C^2$ topology.

The construction of the Markovian return map  in
Section~\ref{ch.rates} 
can be performed in such a way that the following uniformity
condition holds:
 \begin{itemize}
 \item[(u$_0$)]given an integer $N\geq 1$ and  $\e>0$, there
is $\de=\de(\e,N)>0$ such that for $j=1,\dots,N$
 \begin{equation}\label{eq.diferenca}
 \|f-f_0\|_{C^k}<\de\quad \Rightarrow\quad
m\big(\{R_f=j\}\vartriangle\{R_{f_0}=j\}\big)<\e,
 \end{equation}
 where $\vartriangle$ represents the symmetric difference of two
sets.
 \end{itemize}
 This is just by continuity of the inductive construction
for maps in a $C^k$ neighborhood of the original map. In fact, the
construction of the partition on which  the map $R_f$ takes
constant values is based on a finite number of iterations  of the
map $f$. By continuity, we can perform the construction of the
partition in such a way that for some fixed integer $N$ the
Lebesgue measure of $\{R_f=j\}$ varies continuously with the map
$f$ for $j\leq N$. Moreover, the Lebesgue measures of the
auxiliary sets $A_j$ and $B_j$ also vary continuously with the map
$f$ for $j\le N$. Hence, the construction can be carried out with
$R_f$ depending continuously on $f$ as stated in (u$_0$).


\cle\label{le.impu1} Assume  (u$_0$) holds for $f_0$. Suppose
moreover that given any $\e>0$ there are $N\geq 1$ and $\de>0$
 for which
 \begin{equation}\label{eq.soam}
 \|f-f_0\|_{C^k}<\de\quad\Rightarrow\quad
 \big\|\sum_{j=N}^\infty\Chi_{\{R_f>j\}}\big\|_1<\e.
 \end{equation}
 Then uniformity
condition (u$_1$) holds for $f_0$.
 \fle

\dem For the sake of notational simplicity we shall write $R$
instead of $R_f$
 and $R_0$ instead of $R_{f_0}$. We need to show that given $\epsilon>0$ there is  $\de>0$
such that for any $f\in\cf$
 $$
\|f-f_0\|_{C^k}<\de\quad\Rightarrow\quad\|R-R_{0}\|_1<\epsilon.
 $$
 Since
 $$R_0=\sum_{j=0}^{\infty}\Chi_{\{R_{0}>j\}}\qand R=\sum_{j=0}^{\infty}\Chi_{\{R>j\}},$$
 then we have
 \begin{align*}
 \big\|R-R_{0}\big\|_1
  &\leq
 \sum_{j=0}^{N-1}\big\|\Chi_{\{R_{0}>j\}}-
 \Chi_{\{R>j\}}\big\|_1+\big\|\sum_{j=N}^\infty\Chi_{\{R_{0}>j\}}\big\|_1
 + \big\|\sum_{j=N}^\infty\Chi_{\{R>j\}}\big\|_1.
 \end{align*}
 By (u$_0$) and \eqref{eq.soam} all these terms can be made small for $f$ close to $f_0$.
\cqd

Let $\cf$ be a uniform family of non-uniformly expanding maps.
Given $f\in\cf$ we let the expansion time function $\mathcal E^f$
and the recurrence time function $\mathcal R^f$ be defined as in
\eqref{de.exptime} and \eqref{de.slowrec} respectively. The tail
of expansion $\Gamma^f_n$ is also defined for $f\in\cf$ as in
\eqref{tailset} for $n\ge1$.

\cle \label{re.uniform}  Let $\cf$ be a uniform family of $C^k$
($k\ge 2$) non-uniformly maps for which there are $C>0$ and
$\gamma>0$ such that
 $m(\Gamma^f_n)\le Cn^{-\gamma}, $ for all $n\ge 1$ and $f\in\cf$.
 Then the constant $C'$ in Theorem~\ref{t:Markov towers}
 may be taken uniformly in a neighborhood of each $f\in\cf$.
\fle \dem
 As one can easily see from case (B)
in the last part of the previous section, the constant $C'>0$ in Theorem~\ref{t:Markov towers}  depends on the
constant $C>0$. Moreover, from \eqref{e.adelta} and the three possible cases one sees that $C'$  also depends on
some previous constants, namely $\alpha$, $a_1$, $b_1$, $c_1$, $\theta$, $N$ and $R_0$. It is possible to check
that all these constants ultimately depend on the constants $B$, $\beta$, $b$ and $\lambda$ associated to the
non-uniformly expanding map $f$. Naturally they also depend
on the first and second derivatives of $f$. 
We explicit the dependence of the various constants in the table below:

\begin{center}
\begin{tabular}{|c|c|c|}
  \hline
  Constant & Dependence & Reference \\ \hline

 $\sigma$, $\delta$, $\theta$ & $\lambda$ & Proposition~\ref{pr.hyperbolic1}\\

 $\delta_1$  &  $B$, $\beta$, $\sigma$, $\delta$ & Lemma~\ref{leaux}  \\

 $\alpha$ & $\theta$ & \eqref{e.alpha} \\

 $N_0 $ &  $\delta_1$ & \eqref{defN0} \\

 $D_0$, $K_0$ & $\sigma$, $\delta_1$ & Lemma~\ref{le:grow2}\\

 $C_0$ & $B$, $\beta$, $b$, $\sigma$ & Corollary~\ref{co.distortion0}\\

 $C_1$ & $C_0$ & Corollary~\ref{co.distortion}\\

 $\delta_0$ &  $\delta_1$, $\alpha$ & Lemma~\ref{le:grow2}, \eqref{eq.r}, \eqref{e.alpha}\\

 $a_0$ & $\sigma$, $C_1$, $D_0$ & Lemma~\ref{l.flowb} \\

 $b_0$, $c_0$ & $C_1$, $D_0$, $\delta_0$ & Lemma~\ref{l.flowa} \\

 $a_1$ & $a_0$, $b_0$, $c_0$, $\alpha$ & \eqref{e.waa1}, \eqref{e.a1} \\

  $c_1$ &  $\sigma$, $b$, $N_\vare$, $C_0$, $N_0$ & Remark~\ref{re.c1} \\

   $b_1$ & $a_1$, $c_1$ & \eqref{b1}\\

    $c_2$ & $a_1$, $b_0$, $c_0$ & Corollary~\ref{c.omega}\\

 $\vare$ & $K_0$, $N_0$, $\delta_{0}$, $\sigma$ & \eqref{eq.vare} \\

 $N_\vare$ & $\vare$, $\sigma$, $\delta_1$ & Remark~\ref{re.Ne}\\

 $N$ & $N_0$, $N_\vare$ & Proposition~\ref{p.construction}\\

 $R_0$ & $K_0$, $\sigma$, $N_0$, $N$, $\theta$ & Subsection~\ref{sub.exp}, \eqref{e.constants}\\

  \hline
\end{tabular}

\end{center}

\noindent For  better understanding dependencies we use the
convention that no constant depends on a constant from  a line
below. Consequently we have all constants depending on $B$,
$\beta$, $b$ and $\lambda$.
 \cqd

\cpr\label{pr.uniform} Let $\cf$ be a uniform family of $C^k$
($k\ge 2$)
non-uniformly maps for which there are $C>0$ and $\gamma>1$ 
such that
 $m(\Gamma^f_n)\le Cn^{-\gamma}, $ for all $n\ge 1$ and $f\in\cf$.
Then conditions (u$_1$) and (u$_2$) hold for each $f\in \cf$.
 \fpr
 \dem
 Take any $f_0\in\cf$.
If we assume that there are $C>0$ and $\gamma>1$ such that
$m(\Gamma^f_n)\le Cn^{-\gamma}$ for all $n\ge 1$ and all
$f\in\cf$, then by Theorem~\ref{t:Markov towers} there is a
constant $C'>0$ such that $m\{R_f > j\}\le C'n^{-\gamma}$ for all
$n\ge 1$ and all $f\in\cf$,
 as long as $f$ is taken in a sufficiently
 small $C^k$ neighborhood of $f_0$ in $\cf$, say $f\in\cf$ with
 $\|f-f_0\|_{C^k}<\delta$. Actually, as we have observed in Remark~\ref{re.uniform}
 the constant $C'$ may be taken uniformly in a neighborhood of the map $f_0$.
 Thus, given $f\in\cf$ with $\|f-f_0\|_{C^k}<\delta$
 and an integer $N\geq 1$, we have
 $$
 \big\|\sum_{j= N}^\infty\Chi_{\{R_f>j\}}\big\|_1
 \le \sum_{j= N}^\infty m\big(\{R_f > j\}\big)
 \leq \sum_{j= N}^\infty {C'} n^{-\gamma}.
 $$
 Since we are
assuming $\gamma>1$, this last sum can be made arbitrarily small
if we take $N$ large enough. Applying Lemma~\ref{le.impu1} we
obtain uniformity condition (u$_1$).

 For proving that (u$_2$) holds, we  have to show that
the constants $\kappa$ and $K$ in Definition \ref{de.expanding}
may be chosen uniformly for $f$ in a $C^k$ neighborhood of $f_0$
in the uniform family $\cf$. The constant $K$ is given in
Subsection~\ref{ch.rates}.\ref{sub.bd}. As it has been shown
there, it only depends  on $C_0$, $D_0$ and $K_0$. From
Remark~\ref{re.uniform} we see that these constants  may be chosen
uniformly in  $\cf$. On the other hand, the constant $\kappa$
appeared in Subsection~\ref{sub.exp} and depends on  $\sigma$,
$N_0$, $K_0$ and $R_0$, which again may be chosen uniformly in
$\cf$.
 \cqd

As a  consequence of Proposition~\ref{pr.uniform}  and
Theorem~\ref{t.A2} we obtain Theorem~\ref{te.uniform}.


\section{An example}\label{se.examples}

Here we present
robust ($C^1$ open) classes of local diffeomorphisms (with no
critical set) that
 are non-uniformly expanding. Such classes of maps were presented
 in \cite{ABV}, and can be obtained, e.g. through
deformation of a uniformly expanding map by isotopy inside some
small region. In general, these maps are not expanding:
deformation can be made in such way that the new map has periodic
saddles.

 Let $M$ be any compact
manifold supporting some uniformly expanding map $f_0$: there
exists $\sigma_0>1$ such that
$$
\|Df_0(x)v\|>\sigma_0\|v\|   \quad\text{ for every $x\in M$ and
$v\in T_x M$.}
$$
For instance, $M$ could be the $d$-dimensional torus $T^d$. Let
$V\subset M$ be some small compact domain, so that  $f_0\vert V$
is injective. Let $f$ be any  map in a small $C^1$-neighborhood
$\cn$ of $f_0$ so that $\|Df(x)^{-1}\|< \sigma_0$ for every $x$
outside $V$. Assume moreover that the $C^1$-neighborhood
sufficiently small in such a way that:
\begin{enumerate}
\item
 $f$ is {\em volume expanding everywhere}\index{expanding!volume}: there is $\sigma_1>1$ such that
 $$|\det Df(x)| > \sigma_1\quad\text{ for every $x\in M$;}$$
 \item
 $f$ is {\em not too contracting on $V$}: there is some small  $\delta>0$ such that
$$\|Df(x)^{-1}\|< 1+\delta\quad\text{ for every $x\in
V$}.$$
 \end{enumerate}
We are going to show that every map $f$ in such a
$C^1$-neighborhood
$\cn$ of $f_0$ is non-uniformly expanding. 


\cle \label{l.A1} Let $B_1, \ldots, B_p, B_{p+1}=V$ be any partition of $M$ into domains such that $f$ is
injective on $B_j$, for $1\le j \le p+1$. There exists $\theta>0$ (only depending on $f_0$) such that the orbit
of Lebesgue almost every point $x\in M$ spends a fraction $\theta$ of the time in $B_1 \cup \cdots \cup B_p$,
that is, $\#\{0\le j <n: f^j(x)\in B_1 \cup \cdots \cup B_p\}\ge\theta\,n$ for every large $n$. \fle

\dem Let $n$ be fixed. Given a sequence $\underi=(i_0, i_1,
\ldots, i_{n-1})$ in $\{1, \ldots, p+1\}$, we denote
$$
[\underi]= B_{i_0}\cap f^{-1}(B_{i_1})\cap \cdots
                  \cap f^{-n+1}(B_{i_{n-1}}).
$$
Moreover, we define $g(\underi)$ to be the number of values of
$0\le j \le n-1$ for which $i_j \le p$. We begin by noting that,
given any $\theta>0$, the total number of sequences $\underi$ for
which $g(\underi)< \theta \, n$ is bounded by
$$
\sum_{k < \theta \, n} \left(\begin{array}{c}n \\ k\end{array}
\right) p^k 
\le \sum_{k \le \theta \, n} \left(\begin{array}{c}n \\
k\end{array} \right) p^{\theta \, n}
$$
A standard application of Stirling's formula (
gives that the last expression is bounded
by $e^{\gamma n} p^{\theta \, n} 
$, where $\gamma$ depends only on $\theta$ and goes to zero when
$\theta$ goes to zero. On the other hand, since we are assuming
that $f$ is volume expanding everywhere and not too contracting on
$B_{p+1}$, we have $\Leb([\underi])\le \Leb(M) \,
\sigma_1^{-(1-\theta) n}.$ Then the measure of the union $I_n$ of
all the sets $[\underi]$ with $g(\underi)< \theta \, n$ is less
than $ \Leb(M) \sigma_1^{-(1-\theta)n} e^{\gamma n} p^{\theta \,
n}
$ Since  $\sigma_1>1$, we may fix $\theta$ small so that
$e^{\gamma} p^{\theta} 
 < \sigma_1^{1-\theta}$. This means
that the Lebesgue measure of $I_n$ goes to zero exponentially fast
as $n\to\infty$. Thus, by the lemma of Borel-Cantelli, Lebesgue
almost every point $x\in M$ belongs in only finitely many sets
$I_n$. Clearly, any such point $x$ satisfies the conclusion of the
lemma. \cqd

 Let $\theta>0$ be the constant given by Lemma~\ref{l.A1},
and fix $\delta>0$ small enough so that $\sigma_0^{\theta}(1 +
\delta) \le e^{-\lambda}$ for some $\lambda>0$. Let $x$ be any
point satisfying the conclusion of the lemma. Then
$$
\prod_{j=0}^{n-1} \|Df(f^j(x))^{-1}\| \le \sigma_0^{\theta \, n}
(1+\delta)^{(1-\theta)n} \le e^{-\lambda n}
$$
for every large enough $n$. This implies that $x$ satisfies
$$
\limsup_{n\to+\infty}\frac{1}{n}\sum_{j=0}^{n-1}\log\|Df(f^j(x))^{-1}\|
\le -\lambda\,.
$$
and since the conclusion of Lemma~\ref{l.A1} holds Lebesgue almost everywhere we have that $f$ is  a
non-uniformly expanding map.

This shows that any sufficiently small neighborhood of $f$ in the $C^2$ topology constitutes a uniform family of
non-uniformly expanding maps; cf. Definition~\ref{d.uniform}.

\printindex


\end{document}